\documentclass{conm-p-l}

\usepackage{amsmath,enumerate,amsthm,amssymb,amscd}
\input{xy}\xyoption{all}
\CompileMatrices

\newcommand{\gothic}{\mathfrak}
\newcommand{\ra}{\longrightarrow}

\newcommand{\m}{{\gothic{m}}}

\renewcommand{\phi}{\varphi}

\newcommand{\xa}{(x_1^{a_1},\ldots,x_d^{a_d})}

\newcommand{\inc}{\subset}
\newcommand{\Is}{I^{\#}}
\newcommand{\cC}{{\mathcal C}}

\newtheorem{thm}{Theorem}[section]
\newtheorem{lemma}[thm]{Lemma}
\newtheorem{prop}[thm]{Proposition}

\theoremstyle{cor}
\newtheorem{cor}[thm]{Corollary}
\theoremstyle{definition}

\newtheorem{example}[thm]{Example}

\theoremstyle{remark}
\newtheorem{remark}[thm]{Remark}

\numberwithin{equation}{section}

\begin{document}

\title{Minimal Homogenous Liaison and Licci Ideals}

\author{Craig Huneke}

\bigskip

\address{Department of Mathematics \\
        University of Kansas \\
        Lawrence, KS
        66045, USA}

\email{huneke@math.ku.edu}

\urladdr{http://www.math.ku.edu/\textasciitilde huneke}

\author{Juan Migliore}

\address{Department of Mathematics \\
         University of Notre Dame \\
         Notre Dame, IN 46556, USA}

\email{Juan.C.Migliore.1@nd.edu}

\urladdr{http://www.nd.edu/\textasciitilde jmiglior}

\author{Uwe Nagel}

\address{Department of Mathematics \\
         University of Kentucky \\
         Lexington, KY 40506, USA}

\email{uwenagel@ms.uky.edu}

\urladdr{http://www.ms.uky.edu/\textasciitilde uwenagel}

\author{Bernd Ulrich }

\address{Department of Mathematics \\
       Purdue University \\
       West Lafayette, IN
       47907, USA}
\email{ulrich@math.purdue.edu}

\urladdr{http://www.math.purdue.edu/\textasciitilde ulrich}

\subjclass[2000]{13C40, 14M06, 14M10}

\date{\today}

\thanks{The first and last authors were partially supported by the
National Science
Foundation, grants DMS-0244405 and DMS-0200858, respectively. The
second and third authors were partially supported by the National
Security Agency, grants H98230-07-1-0036 and H98230-07-1-0065,
respectively.}

\begin{abstract}
We study the linkage classes of homogeneous ideals in polynomial
rings. An ideal is said to be homogeneously licci if it can be
linked to a complete intersection using only homogeneous
regular sequences at each step. We ask a natural question:
 if $I$ is homogeneously licci, then
can it be linked to a complete intersection  by linking using
regular sequences  of forms  of smallest  possible  degree at each
step (we call such  ideals minimally homogeneously licci)? In this
paper we answer this question in the negative. In particular, for
every $n\geq 28$ we construct a set of $n$ points in $\mathbb P^3$
which are homogeneously licci, but not minimally homogeneously
licci. Moreover, we prove that one cannot distinguish between the
classes of homogeneously licci and non-licci ideals based only on
their Hilbert functions, nor distinguish between homogeneously licci
and minimally homogeneously licci ideals based solely on the graded
Betti numbers. Finally, by taking hypersurface sections, we show
that the natural question has a negative answer whenever the
height of the ideal is at least three.
\end{abstract}

\maketitle




\section{Introduction}

\medskip

Let $R$ be a commutative Noetherian ring, and let $I$ and $J$ be two
proper ideals of height $c$ in $R$. These ideals are said to be \it
directly linked \rm if there exists a regular sequence
$f_1,\ldots,f_c$ such that $(f_1,\ldots,f_c):I = J$ and
$(f_1,\ldots,f_c):J = I$. We say $I$ and $J$ are in the same \it
linkage class \rm (or \it liaison class\rm) if there exists a
sequence of ideals $I = I_0, \ldots, I_n = J$ such that $I_j$ is
directly linked to $I_{j+1}$ for $0\leq j\leq n-1$, the case $n=2$
being referred to as {\it double linkage}. Such a sequence of links
connecting $I$ and $J$ is far from unique. We call the ideal $I$ \it
licci \rm if $I$ is in the {\bf li}nkage {\bf c}lass of a {\bf
c}omplete {\bf i}ntersection, i.e., of an ideal generated by a
regular sequence.

Let $I$ be a homogeneous ideal of $S = k[x_1,\ldots,x_d]$, a polynomial
ring over a field with the standard grading. In considering the linkage
class of $I$ one must distinguish between allowing only homogeneous links,
i.e., links where all the regular sequences are homogeneous, and
non-homogeneous
linkage. The corresponding linkage classes are not a priori the same.
In trying to understand homogeneous linkage, an obviously critical question
to address is to understand the sequence of links by which a
homogeneously licci ideal is
in the liaison class of a complete intersection.
In the local case,
in \cite[2.5]{HUduke} it was shown
that if an ideal $I$ in a Gorenstein local ring with infinite
residue field  is licci, then  one  can  pick  regular  sequences
consisting of general combinations of minimal generators  at every  step
to  reach a complete  intersection. A corresponding question for homogeneous
linkage would ask
if $I$ is homogeneously licci, then
can it be linked to a complete intersection  by linking using  regular
sequences  of forms  of smallest  possible  degree at each step?
In this paper we answer this question in the negative.

Since homogeneous linkage by a regular sequence of minimal possible
degree will be the chief topic of this paper, we give it a name; we
say that $I$ is {\it directly minimally homogeneously linked} to an
ideal $J$ if $J = L:I$ and $I = L:J$ where $L$ is generated by a
regular sequence in $I$ of smallest possible degree. Note that
contrary to usual linkage, this is not necessarily a symmetric
relation!  We say that $I$ is {\it minimally homogeneously licci} if
$I$ is homogeneously licci and all  the  links in a linkage sequence
from $I$ to a complete intersection can be taken to be direct
minimal homogeneous links. (Of course, this does not mean \it all
\rm paths from $I$ to a complete intersection are minimal linkages.)

Intuitively, it seems that all homogeneously licci ideals should
be minimally homogeneously licci, as is the case in height $2$.
However, in  \cite[pg. 61-62]{GHMS}
an example was announced of a homogeneously licci ideal which is not
minimally
homogeneously licci. An $\m$-primary homogeneous ideal $I$ with
h-vector $(1,3,6,8,7,6,2)$ was constructed which is homogeneously licci,
and it
was claimed that it is not minimally homogeneously licci. However,
no proof was given for this assertion. We began the present paper in an
attempt to understand this
example. In Section 2 we will give examples of three ideals,
each with h-vector $(1,3,6,8,7,6,2)$, such that one of them is not
licci, one is minimally homogeneously licci, and the last  is
homogeneously licci, but not minimally homogeneously licci. Although
none of these examples is the one from \cite{GHMS} (which was level,
while ours are not), they show
that no argument based only on the Hilbert function
can distinguish such ideals. In the last two examples,
even the graded Betti numbers are the same. This proves that there
is no way to distinguish homogeneously licci and minimally homogeneously
licci ideals by only using their graded Betti
numbers. At  this point, we do not know  whether there exist two
homogeneous ideals with the same graded Betti numbers, one
homogeneously licci and the other not.

In our last section we go much farther: we prove that for any $n\geq
28$ and $d\geq 3$ there exist reduced subschemes of degree $n$ and
codimension three in $\mathbb P^d_k$, for $k$ an infinite
field, that are homogeneously licci, but not minimally
homogeneously licci. Moreover, these subschemes can be taken
to be finite unions of reduced linear subspaces. In particular,
for any $n \geq 28$ there exist
sets of $n$ $k$-rational points in $\mathbb P^3_k$ that are
homogeneously licci, but are not minimally homogeneously licci. The
examples given are quite explicit, and suggest that this phenomenon
is in fact very common. Finally, by taking hypersurface sections, we
show that in any codimension $\geq 3$ there are reduced subschemes
that are homogeneously licci, but not minimally homogeneously licci.
Related results for Gorenstein liaison recently have been given in
\cite{HSS}. On the other hand, it is still conceivable that every
homogeneously licci ideal can be linked to a complete intersection
by using only homogeneous regular sequences that are part of a
minimal homogeneous generating set of the respective ideals.

For unexplained terminology we refer to the book of Eisenbud,
\cite{E}, while for additional information on liaison we refer to
\cite{Mi}.

We would like to thank the referee of this article for a careful
reading and many insightful comments.

\bigskip
\smallskip

\section{Three $\m$-Primary Examples}

\medskip

In this section we give three examples of $\m$-primary homogeneous
ideals, all with Hilbert function $(1,3,6,8,7,6,2)$, which exhibit
different behavior in terms of their linkage classes. One is not
licci at all, one is minimally homogeneously licci, and the third is
homogeneously licci, but not minimally homogeneously licci. The
third example takes considerably more work.

We begin by establishing some notation. We will always write $S =
k[x_1,\ldots,x_d]$ for a polynomial ring over a field $k$ and $\m$
for its homogeneous maximal ideal $(x_1,\ldots,x_d)$. If $I$ is a
homogeneous
ideal, we let $I_{\leq j}$ denote the ideal generated by all forms of
degree at most $j$ in $I$.

One of our examples is monomial. By a
monomial in $S$ we mean an element of the form $x_1^{a_1}\cdots
x_d^{a_d}$.
A {\it monomial ideal} is an ideal generated by monomials. Every
$\m$-primary monomial ideal $I$ can be written uniquely in {\it
standard form} \, $I = \xa + I^{\#}$, where $\Is$ is generated by
monomials that together with $\{x_1^{a_1},\ldots,x_d^{a_d}\}$
generate $I$ minimally.

We will use the following lemma of  \cite[2.5]{HU}:

\begin{lemma}\label{cancel} Let $S = k[x_1,\ldots,x_d]$
be a polynomial ring over a field $k$ and let $I$ be an
$\m$-primary monomial ideal. If $\Is = x_1^{b_1}\cdots x_d^{b_d}K$
where $K\neq S$ is a monomial ideal,
then the ideal $I' = (x_1^{a_1-b_1},\ldots, x_d^{a_d-b_d}) + K$ is
obtained from $I$ by a double link defined by the monomial regular
sequences $x_1^{a_1},\ldots,x_d^{a_d}$ and
$x_1^{a_1-b_1},\ldots,x_d^{a_d-b_d}$.
\end{lemma}

\smallskip

\begin{example}\label{examp1}{\rm Let  $I =
(z^3, xyz, x^3y, x^4, y^6, y^5z, xy^5)\subset S=k[x,y,z]$. Then $I$
is $\m$-primary and the
Hilbert function of $S/I$ is $(1,3,6,8,7,6,2)$. We claim that
$I$ is not licci. We compute $I^{\#} = (xyz, x^3y, y^5z, xy^5)$.
This  ideal has height one with greatest  common divisor $y$. By
Lemma~\ref{cancel}, $I$ is  doubly  linked to $I' = (xz, x^3, y^4z,
xy^4, z^3, y^5)$, and $(I')^{\#} = (xz, y^4z, xy^4)$, which has
height  two. Now \cite[2.4]{HU} shows that $I$ is
not  licci  in  any sense. We note that the minimal homogeneous free
resolution of  $I$ is
$$0\longrightarrow S(-8)^2\oplus S(-9)^2\longrightarrow S(-5)^3\oplus
S(-7)^5\oplus S(-8)^2
\longrightarrow S(-3)^2\oplus S(-4)^2\oplus S(-6)^3.$$
}
\end{example}

\smallskip

\begin{example}\label{examp2}{\rm Let $I = (z^3, xyz, x^3y, x^4, y^6,
y^5z+x^2y^4, xy^5)
\subset S=k[x,y,z]$. Then $I$ is an $\m$-primary ideal with exactly the same
Hilbert function as the
ideal in Example~\ref{examp1}. However,  $I$ is  licci,  and  it  is
even minimally homogeneously licci. We  prove this by  simply
constructing a sequence of minimal homogeneous links  to a complete
intersection. Let $L = (z^3, x^4, y^6)$. Then $L$ is a complete
intersection inside  $I$ of degrees $3,4,6$, which are the minimal
possible ones since  $I_{\leq 3}$  has  height  one, and $I_{\leq 5}$
has  height two. A calculation  gives that $I_1: = L:I = (z^3,
x^3y, x^4, x^3z-xyz^2, y^5)$. The minimal possible degrees of a
homogeneous regular sequence in $I_1$ are $3,4,5$. Let $L_1 = (z^3, x^4,
y^5)$
and set $I_2: = L_1:I_1 = (xz, x^3, z^3, xy^4, y^5, x^2y^3+y^4z)$.
The  minimal possible  degrees of a homogeneous regular sequence in
$I_2$ are
$2,3,5$. Set $L_2 = (xz, x^3+z^3, y^5)$. A calculation  gives that
$I_3: = L_2:I_2 = (xz, yz, xy+z^2, x^3, y^5)$. The minimal possible
degrees of a homogeneous regular sequence in $I_3$ are $2,2,5$. Let $L_3
= (yz,
xy+z^2, y^5+x^5)$ and set $I_4: = L_3:I_3 =
 (yz, xy+z^2, y^2, x^4, x^2z^2)$. The  minimal
possible  degrees of a homogeneous regular sequence in $I_4$ are
$2,2,4$. Let
$L_4 = (xy+z^2, y^2, x^4)$ and set $ I_5 = L_4:I_4 = (y, xz, z^2, x^4)$.
The  minimal
possible  degrees of a homogeneous regular sequence in $I_5$ are
$1,2,4$. Set
$L_5 = (y,  z^2, x^4)$. Then $L_5:I_5 = (y,z,x^3)$ is a complete
intersection.}

We note that the minimal homogeneous free resolution of $I$ is
$$0\longrightarrow S(-8)\oplus S(-9)^2\longrightarrow S(-5)^3\oplus
S(-7)^5\oplus S(-8)
\longrightarrow S(-3)^2\oplus S(-4)^2\oplus S(-6)^3 .$$
\end{example}

\bigskip
\smallskip

The first example was easy to verify using the work of \cite{HU}.
The second example was easily verified by a straightforward
calculation. However, the ease of proving the second example belies
the difficulty in finding it. We needed to find a licci ideal with
the given Hilbert function where a certain Koszul relation on the
generators is minimal. The latter feature is crucial in the
construction, as it causes Betti numbers to decrease under double
linkage.

Our last example is more
delicate, requiring an understanding of how the heights of ideals
generated by
forms of low degree change in the linkage class of the ideal.

\smallskip

\begin{thm}\label{mainthm} Let $S=k[x,y,z]$, where $k$ is an infinite
field, and consider
the ideal $I = (x^2y+y^3-yz^2-z^3, xy^2-xz^2, x^3z-xyz^2,
y^2z^2-z^4, x^6, z^6, xz^5)$. Then $I$ is an $\m$-primary ideal with
minimal homogeneous free resolution
$$0\longrightarrow S(-8)\oplus S(-9)^2\longrightarrow S(-5)^3\oplus
S(-7)^5\oplus S(-8)
\longrightarrow S(-3)^2\oplus S(-4)^2\oplus S(-6)^3,$$
the Hilbert function  of $S/I$ is $(1,3,6,8,7,6,2)$, and
$I$ is homogeneously licci, but not minimally homogeneously licci.
\end{thm}

\begin{proof} We first prove that $I$ is homogeneously
licci by direct computation. Let $L$ be generated by the
homogeneous regular sequence $x^2y+y^3-yz^2-z^3, y^2z^2-z^4,  x^6$.  A
computation shows that
$I_1: = L:I = (x^2y+y^3-yz^2-z^3, z^4, y^2z^2, xz^3, x^5)$.  We now let
$L_1$ be the
ideal generated by the homogeneous regular sequence $x^2y+y^3-yz^2-z^3,
z^4, x^5$.
One obtains the link $I_2: = L_1:I_1 = (z^2, x^2y+y^3, x^2z+y^2z, y^5,
xy^3z, x^5)$.
Let $L_2$ be the ideal in $I_2$ generated by the homogeneous regular
sequence $z^2, x^2y+y^3, x^5$.
One can compute $I_3: = L_2:I_2 = (xz, y^2, z^2, x^2y, x^5).$ We let
$L_3$ be the ideal
in $I_3$ generated by the regular sequence $y^2, z^2, x^5$, and
calculate $I_4: = L_3:I_3 =
(y^2, yz, z^2, x^4,x^3z)$. Let $L_4$ be generated by the regular sequence
$y^2, z^2, x^4$, and $I_5: = L_4:I_4 = (z, xy, y^2, x^4)$.
Set $L_5 = (z,y^2,x^4)$. Then $I_6: = L_5:I_5 = (z, y, x^3)$ is
generated by a regular sequence,
so that $I$ is licci. The very first link, however, is not minimal; it
uses a homogeneous regular
sequence of degrees $3,4,6$, whereas the minimal degrees of a
homogeneous regular sequence
inside $I$ are $3,3,6$.

 The claims concerning the Hilbert function  of $I$  and
the  resolution of  $I$  are  easy  to  check on  any computer
algebra program. The remainder of the proof is to show that $I$ is
not minimally homogeneously licci.

We consider any $\m$-primary homogeneous ideal $I'$ in
$S=k[x,y,z]$ satisfying condition
$(\star)$; by this latter condition  we mean that
$I'$ contains a homogeneous regular sequence of degrees $3,3,6$ and
has a (not necessarily minimal) homogeneous free resolution of the
form
$$0\longrightarrow S(-8)\oplus S(-9)^2\longrightarrow S(-5)^3\oplus
S(-7)^5\oplus S(-8)
\longrightarrow S(-3)^2\oplus S(-4)^2\oplus S(-6)^3.$$

The ideal $I$ satisfies $(\star)$. Also notice that any ideal
satisfying $(\star)$ requires $7$ generators, hence cannot be a
complete intersection or directly homogeneously linked to a complete
intersection. We are going to prove that any minimal homogeneous
double link reproduces condition $(\star)$. Therefore no ideal
satisfying $(\star)$ can be minimally homogeneously licci.

Thus assume that $I'$ satisfies $(\star)$. We first prove that $H$,
the subideal of $I'$ generated by all forms of degree at most $5$,
has height at most $2$. Indeed if $H$ has height three, then $H$
contains an ideal $L$ generated by a homogeneous regular sequence of
degrees $3,3,4$.
The $k$-dimension of $L_7$ is $35$. On the other hand the resolution
of $I'$ shows that the $k$-dimension of $H_7$ is at most
$2\dim_k(S(-3)_7) + 2\dim_k(S(-4)_7)- 3\dim_k(S(-5)_7) = 32$. This
contradiction proves that $H$ has height at most two. Thus the
smallest degrees of a homogeneous regular sequence contained in $I'$
is $3,3,6$.

Let $J$ be a link of $I'$ using a homogeneous regular sequence of degrees
$3,3,6$. We obtain its resolution from Ferrand's mapping cone construction
(cf. \cite{PS}). Since $3,3,6$ are the minimal possible degrees of a
homogeneous
regular sequence inside $I'$, any such regular sequence is part of a
minimal homogeneous generating set of $I'$; in particular, we can
assume that this regular sequence appears in the generating set
given by the above resolution of $I'$. This accounts for splitting,
in the mapping cone, of $S(-3)^2\oplus S(-6)$. Thus we obtain the following
(not necessarily minimal) homogeneous free resolution of $J$:
$$0\longrightarrow S(-6)^2\oplus S(-8)^2\longrightarrow S(-4)\oplus
S(-5)^5\oplus
S(-7)^3 \longrightarrow S(-3)^4\oplus S(-4)\oplus S(-6).$$ If
$3,3,6$ are again the minimal degrees of a homogeneous regular
sequence contained in $J$, then linking once more we reproduce an
ideal satisfying condition $(\star)$ as asserted. Otherwise $J$
contains a homogeneous regular sequence of degrees $3,3,4$.
Linking $J$ with respect to such a complete intersection would
lead to this minimal homogeneous free resolution of the linked
ideal $K$,
$$...\rightarrow S(-3)^{a+1}\oplus S(-5)^5\oplus
S(-6)^b \rightarrow S(-2)^2\oplus S(-3)^{a}\oplus
S(-4)^3\rightarrow K\rightarrow 0,$$ where $a\geq 0$ and $b \geq 0$.
Since the two linearly independent quadrics in $K$
must have exactly one generating syzygy, it follows that $a+1=1$,
hence $a=0$. Thus
$K_{\leq 3}$ is generated by these two quadrics. Furthermore, having
a nontrivial cubic relation, the two quadrics are forced to have a
linear factor in common. We conclude that $K_{\leq 3}$ has height
one and hence cannot contain two cubics forming a regular sequence.
Thus there exists no homogeneous regular sequence of degrees $3,3,4$
inside $J$. This finishes the proof.
\end{proof}

\smallskip

One crucial difference between the examples in (2.3) and (2.4) is
that in the first the cubics in the ideal have height one, while in
the second they have height two.

\bigskip
\medskip

\section{Minimal Homogeneous Linkage in $\mathbb P^d$}
\medskip

In this section we continue our discussion of homogeneously licci,
but not minimally homogeneously licci ideals, by constructing for
any $n\geq 28$ a reduced subscheme of $\mathbb P^d$ of degree $n$
whose defining ideal is homogeneously licci, but not minimally
homogeneously licci. We begin with a general remark.

\begin{remark}\label{socle} Let $S = k[x_1,\ldots,x_d]$ be a polynomial
ring over
a field $k$, and let $I$ be a homogeneous ideal in $S$ of
height $c$ such that $S/I$ is Cohen-Macaulay. Let
$$0\rightarrow \oplus_i S(-i)^{b_{c,i}}\rightarrow \oplus_i
S(-i)^{b_{c-1,i}}\rightarrow
\ \ldots \ \rightarrow \oplus_i S(-i)^{b_{1,i}}\rightarrow
I\rightarrow 0$$ be a minimal homogeneous free resolution of $I$.
Suppose that $I$ contains an ideal generated by a homogeneous
regular sequence, say $f_1,\ldots,f_c$, of degrees $d_1,\ldots,d_c$,
respectively. Then
$$ d_1+\cdots + d_c \geq \text{max}\{i|\,\,b_{c,i}\ne 0\}.$$
Furthermore, if $I$ is not equal to $(f_1,\ldots, f_c)$, then the inequality
above is strict.
\begin{proof}
Without loss of generality we may assume that $k$ is infinite.
Choose a sequence of linear forms $l_1,\ldots,l_{d-c}$ which form a
regular sequence on both $S/I$ and $S/(f_1,\ldots,f_c)$. The image
of $I$ in $T = S/(l_1,\ldots,l_{d-c})$ has the same graded Betti numbers as
$I$, and the images of $f_1, \ldots , f_c$
in $T$ form a regular sequence. Since both $T/IT$ and
$T/(f_1,\ldots,f_c)T$ have finite length, it follows from
\cite[Exercises 20.18 and 20.19]{E} that the regularity is given by the
formulas
$$
\text{reg}(T/(f_1,\ldots,f_c)T) = d_1+\cdots + d_c - c =
\text{max}\{n|\, ((T/(f_1,\ldots,f_c)T)_n\ne 0\}$$ and
$$\text{reg}(T/IT) = \text{max}_i\{i|\,\,b_{c,i}\ne 0\}-c =
\text{max}\{n|\, (T/IT)_n\ne 0\}.$$

Clearly the largest $n$ for which $(T/IT)_n\ne 0$ must be at most
the largest $n$ for which $(T/(f_1,\ldots,f_c)T)_n\ne 0$. Thus,
$d_1+\cdots + d_c \geq \text{max}\{i|\,\,b_{c,i}\ne 0\}$. Moreover,
if $I\ne (f_1,\ldots,f_c)$, then this remains true after passing to
$T$, since $l_1,\ldots,l_{d-c}$ form a regular sequence modulo $I$.
In this case, since $T/(f_1,\ldots,f_c)T$ is Gorenstein,
$(f_1,\ldots, f_c)T:\m T \subset IT$, and therefore the largest $n$ for
which $(T/IT)_n\ne 0$ must be strictly smaller than the largest $n$
for which $(T/(f_1,\ldots,f_c)T)_n\ne 0$. Hence $d_1+\cdots + d_c >
\text{max}_i\{i|\,\,b_{c,i}\ne 0\}$. \end{proof}
\end{remark}

\medskip

The next theorem provides a large class of examples of subschemes
in $\mathbb P^d$ which are homogeneously licci, but not minimally
homogeneously licci.
The subschemes we obtain have codimension $3$, which is the smallest
possible codimension for such examples.

\medskip

\begin{thm} \label{two bdl} Let $S = k[x_0,\ldots,x_d]$ with $d\geq 3$,
where $k$ is an infinite field.
Let $4 \leq a_1+3\leq  a_2 < a_3 < a_4$ be integers so that $a_1
\neq 2$ and $a_2+a_3\leq {\rm min}\{2,a_1 \} +a_4$.
 Choose $F_1, F_4\in S$ homogeneous elements of degrees $a_1, a_4$,
 respectively, and a
 linear form $L_1$ such that $L_1, F_1, F_4$ is a regular sequence.
 Define $I_{1} = (L_1, F_1, F_4)$.
Choose $F_2, F_3
 \in I_{1}$ homogeneous elements of degrees $a_2, a_3$, respectively, so
that $F_2, F_3$ and $L_1, F_2$
 form regular sequences. $($This is possible since $(L_1, F_1)$ has
height $2$. Necessarily
$(F_2,F_3)\inc (L_1, F_1)$.$)$ Let $L_2$ be a linear form such that
$L_2, F_2, F_3$ is a regular sequence.
 Define $I = L_2 \cdot I_{1} + (F_2, F_3)$. Then $Z = V(I)$ is
homogeneously licci, but not
 minimally homogeneously licci.
\end{thm}

\begin{proof}
We begin by verifying the assertion that $I$ is licci. In fact it is
doubly  linked
to $I_1$, a complete intersection. To see this, choose any
homogeneous $G\in I_1$ such that $G, F_2, F_3$ form a regular
sequence. Since the  linear form $L_2$ is regular modulo
$(F_2,F_3)$ one sees that $J: = (G, F_2, F_3):I_1 =
(GL_2, F_2, F_3):I$. Moreover, as the calculation of a free resolution
of $I$ done
below proves, $S/I$ is Cohen-Macaulay, hence $I$ is unmixed and $(GL_2,
F_2, F_3):J = I$.

Next we find a free resolution of $I$. The definition of $I$
immediately gives that $I:L_2=I_1$, again as the linear form $L_2$
is regular modulo $(F_2,F_3)$.
Thus there is an exact sequence
$$0\ra S/(L_1, F_1, F_4)(-1)\ra S/I\ra S/(L_2, F_2, F_3)\ra 0 \, .$$
Using a Koszul resolution of the first and the last module in this
sequence, applying the horseshoe lemma, and splitting off a summand we
obtain a (not
necessarily minimal) homogeneous free resolution of $I$:

\[
0 \rightarrow
\begin{array}{c}
S(-a_1 -a_4-2) \\
\oplus \\
S(-a_2 -a_3 -1)
\end{array}
\rightarrow
\begin{array}{c}
S(-a_1 -2) \\
\oplus \\
S(-a_4 -2) \\
\oplus \\
S(-a_1 -a_4 -1) \\
\oplus \\
S(-a_2 -a_3) \\
\oplus \\
S(-a_2 -1) \\
\oplus \\
S(-a_3 -1)
\end{array}
\rightarrow
\begin{array}{c}
S(-2) \\
\oplus \\
S(-a_1 -1) \\
\oplus \\
S(-a_4-1) \\
\oplus \\
S(-a_2) \\
\oplus \\
S(-a_3)
\end{array}
\rightarrow I \rightarrow 0 \, .
\]
Observe that the degrees of the generators are ordered by $2 \leq a_1 +1
< a_2 < a_3 < a_4+1$.

Let $I'$ be an ideal satisfying the following three conditions,
which we denote by $(\star)$: $I'$ has height three, contains a
homogeneous regular sequence of degrees $2,a_2,a_4+1$, and has a
(not necessarily minimal) homogeneous free resolution of the form
\[
0 \rightarrow
\begin{array}{c}
S(-a_1 -a_4-2) \\
\oplus \\
S(-a_2 -a_3 -1)
\end{array}
\rightarrow
\begin{array}{c}
S(-a_1 -2) \\
\oplus \\
S(-a_4 -2) \\
\oplus \\
S(-a_1 -a_4 -1) \\
\oplus \\
S(-a_2 -a_3) \\
\oplus \\
S(-a_2 -1) \\
\oplus \\
S(-a_3 -1)
\end{array}
\rightarrow
\begin{array}{c}
S(-2) \\
\oplus \\
S(-a_1 -1) \\
\oplus \\
S(-a_4-1) \\
\oplus \\
S(-a_2) \\
\oplus \\
S(-a_3)
\end{array}
\rightarrow I' \rightarrow 0 \, .
\]

The ideal $I$ satisfies $(\star)$ since $L_1L_2, F_2$ form a homogeneous
regular sequence in $I$
of degrees $2, a_2$ and since $I$ is generated in degrees
at most $a_4+1$.

We will show that no ideal $I'$ satisfying $(\star)$ is a complete
intersection, nor is any direct minimal homogeneous link of $I'$.
We will also prove that minimal homogeneous double linkage
reproduces condition $(\star)$. Therefore no ideal satisfying
$(\star)$ can be minimally homogeneously licci.

We first argue that no ideal $I'$ satisfying $(\star)$ can be
generated by a regular sequence. Indeed, the two smallest degrees of
minimal generators of $I'$ are $2, a_1+1$, and there is a minimal
syzygy in degree $a_1+2$. The fact that the syzygy in degree $a_1+2$
is minimal follows from our assumptions, most notably the inequality
$a_1+3\leq a_2 $, which imply that there is no cancelation between
this syzygy and a generator or a second syzygy of $I'$.

Next we claim that the minimal possible degrees of a homogeneous regular
sequence in $I'$  are $2, a_2, a_4 +1$. By assumption $I'$ contains a
homogeneous regular sequence of these degrees. We need to
show that $I'$ cannot have a smaller regular sequence.
First suppose that the ideal $I'_{\leq a_1+1}$ has
height two. As noted above, there is a minimal first syzygy of
degree $a_1+2$, which can only come from the two elements in $I'$ of
degrees $2$ and $a_1+1$. But then these elements cannot form a regular
sequence, proving the claim. It remains to show that $I'$ does not contain
a homogeneous regular sequence of degrees $2, a_2, a_3$. Suppose that
such a regular
sequence exists. By Remark \ref{socle}, it follows that
$$2+a_2+a_3 > \text{max}\{a_1+a_4+2, a_2+a_3+1\} = a_1+a_4+2$$
(recall that by assumption $a_2+a_3\leq a_1+a_4$). But then
$a_2+a_3 > a_1+a_4$, a contradiction.

We now compute a homogeneous free resolution of any link $J$ of $I'$
with respect to a homogenous regular sequence of degrees $2, a_2,a_4
+1$. Since such a sequence is part of a minimal homogeneous generating
set of $I'$, copies of $S(-2)$, $S(-a_2)$ and $S(-a_4-1)$ split off in
the mapping cone construction (cf. \cite{PS}). Thus we obtain the following
homogeneous free resolution of $J$:

\[
0 \rightarrow
\begin{array}{c}
S(a_3 - a_2 - a_4 -3) \\
\oplus \\
S(a_1 - a_2 - a_4 - 2)
\end{array}
\rightarrow
\begin{array}{c}
S(a_1 - a_2 - a_4-1) \\
\oplus \\
S(-a_2 -1) \\
\oplus \\
S(a_1 - a_2 -2) \\
\oplus \\
S(a_3 - a_4 -3) \\
\oplus \\
S(-a_4 -2) \\
\oplus \\
S(a_3 - a_2 - a_4 - 2)
\end{array}
\rightarrow
\begin{array}{c}
S(a_1 - a_2 - 1) \\
\oplus \\
S(a_3 - a_4 - 2) \\
\oplus \\
S(-a_4 -1) \\
\oplus \\
S(-a_2) \\
\oplus \\
S(-2)
\end{array}
\rightarrow J \rightarrow 0 \, .
\]

We first argue that $J$ cannot be a complete intersection. Indeed,
the summands $S(a_1-a_2-a_4-2)$ and $S(a_1-a_2-a_4-1)$ in the last
and next to last step of the above resolution are still present
after passing to a minimal homogeneous free resolution, due to our
numerical conditions. Hence the
minimal homogeneous free resolution of $S/J$ cannot be symmetric as
$J$ does not contain a linear form. It follows that $S/J$ is not
even Gorenstein.

Notice that the degrees of generators are $a_2 - a_1 +1 $, $a_4 -
a_3+2$, $a_4 +1$, $a_2$, $2$, and that $2 < a_2 - a_1 + 1 \leq a_2
\leq a_4-a_3+2 < a_4 +1$ according to our numerical assumptions. By
construction, $J$ contains a homogeneous regular sequence of degrees
$2, a_2, a_4+1$. We claim that these are the smallest possible
degrees of a homogeneous regular sequence in $J$. Suppose the
contrary. In this case either $a_1
 > 1$ and $J$ contains a homogeneous regular sequence of degrees $2,
a_2-a_1+1$,
or else $J$ contains a homogeneous regular sequence of degrees $2,
a_2, a_4 - a_3 + 2$.

To rule out the first possibility notice that the two generators of
degrees $2$ and $a_2-a_1+1$ are the only generators of degrees $<
a_2-a_1+2$ since $a_2-a_1+1 < a_2$. Furthermore, the syzygy of
degree $a_2-a_1+2$ cannot cancel against any second syzygy or
against any generator of $J$, because $a_1 >2$. Thus this syzygy is
a minimal syzygy, and hence gives a nontrivial relation of degree
$a_2-a_1+2$ among the two generators of degrees $2$ and $a_2-a_1+1$. It
follows that these generators could not form a regular sequence.

To rule out the second possibility from above, suppose that
$J$ contains a homogeneous regular sequence of degrees $2, a_2, a_4 -
a_3 +2$.
Now Remark~\ref{socle} implies that
$$2+a_2+a_4-a_3+2 > \text{max}\{a_2+a_4-a_3+3, a_2+a_4-a_1+2\} =
a_2+a_4-a_1+2.$$
But this is impossible by our numerical conditions.

Thus, when we perform a minimal homogeneous link for $J$, the residual again
satisfies condition $(\star)$, proving the theorem. \end{proof}

\smallskip

\begin{cor} \label{p leq 10} Let $k$ be an infinite field.
For any $n \geq 28$ and $d\geq 3$ there exist reduced subschemes of
degree $n$ and codimension three in $\mathbb P^d_k$ that are
homogeneously licci, but not minimally homogeneously licci.
Moreover, the subschemes can be chosen to be a finite union of
reduced linear subspaces. In particular, for any $n \geq 28$ there exists a
set of $n$ $k$-rational points in $\mathbb P^3_k$ that is homogeneously
licci,
but not minimally homogeneously licci.
\end{cor}

\begin{proof}
The last statement follows at once from the case $d=3$.

We will use Theorem~\ref{two bdl}  with $a_1=1$. Let $l_2,\ldots
,l_{a_2}, m_2, \ldots ,m_{a_3}, n_1, \ldots, n_{a_4}$, $L_1$, $L_2$
be linear forms
in $S$, which together with $x_0$ have the
property that any four are linearly independent.
Such linear forms exist because $d \geq 3$ and $k$ is infinite.
Setting
$$F_1=x_0, F_2= x_0\cdot l_{2}\cdots l_{a_2}, F_3= L_1\cdot m_2\cdots
m_{a_3}, F_4= n_1\cdots n_{a_4}$$
we obtain homogeneous elements in $S$ that satisfy the assumptions of
Theorem~\ref{two bdl}. Moreover,  $L_2$ is regular modulo $I_1$.
Thus the definition of $I$ shows that
$$
I = (L_2,F_2,F_3) \cap I_1 = (L_2,F_2,F_3) \cap (L_1,F_1,F_4),
$$
and each primary component of this ideal is generated by linear
forms, as can be seen from our linear independence condition.
Therefore $V(I)$ is a finite union of reduced linear subspaces,
of degree $a_2a_3+a_4$.
Moreover, $V(I)$ is homogeneously licci, but not minimally homogeneously
licci if the numerical
conditions of the theorem are satisfied. Choosing $a_1 = 1, a_2 = 4,
a_3 = 5, a_4 = 8$ gives a subscheme of degree $28$, and simply
increasing the value of $a_4$ gives all other possible $n$.
\end{proof}

\bigskip

We were able to prove a stronger result by replacing the lower bound
$28$ for the number of points by the value $10$, but the proof was
more complicated and we opted to present the simpler proof with a
slightly worse lower bound. We do not know what the best lower bound
is.

Using the above result we are now going to show that in any
codimension $\geq 3$, there exist homogeneously licci subschemes
that are not minimally homogeneously licci. We continue to assume
that $k$ is infinite. If $C \subset I$ are ideals
we call the colon ideal $C:I$ the {\it residual} of $I$ with respect to $C$.

\smallskip

\begin{prop} \label{hypersurf sect}
Consider a class $\cC$ of proper homogeneous ideals of height $c$ in $S =
k[x_0,\ldots,x_d]$ that is closed under taking residuals with respect to
homogeneous complete intersections of type $(d_1,\ldots,d_c)$. Assume further
that $(d_1,\ldots,d_c)$ is the minimal
possible type of a homogeneous complete intersection of height $c$ contained in
any ideal in $\cC$.

Let $I$ be a fixed homogeneously licci ideal in $\cC$ and let $F \in
S$ be a general form of sufficiently large degree $D$ that is
regular modulo $I$. Then the ideal $(I, F)$ is homogeneously licci,
but not minimally homogeneously licci.

Moreover, consider the class $\cC'$ of ideals of the form $(I',F')$,
where $I' \in \cC$ and $F' \in S$ is a form of degree $D$ that is
regular modulo $I'$. Then $\cC'$ is closed under taking residuals
with respect to homogeneous complete intersections of type $(d_1,
\ldots, d_c, D)$, and $(d_1, \ldots, d_c, D)$ is the minimal
possible type of a homogeneous complete intersection of height $c+1$
contained in any ideal in $\cC'$.
\end{prop}

\begin{proof}
Since $I$ is homogeneously licci, there is a finite sequence of
homogeneous links taking $I$ to a complete intersection.  Since $F$
is general and $k$ is infinite, $F$ is not only regular modulo $I$
but also modulo all these intermediate ideals.  It then follows
that $(I, F)$ is homogeneously licci (see
\cite[5.2.17]{Mi}).

The asserted properties of $\cC'$ imply that a complete intersection of height $c+1$
in $\cC'$ must be of type $(d_1,\ldots,d_c,D)$, hence its residual with respect
to itself is not a proper ideal and could not belong to $\cC'$. It follows that no ideal
in $\cC'$ can be a complete intersection. Now again the asserted properties of $\cC'$
give that $(I,F)$ is not minimally homogeneously licci.

It remains to prove the two claims about $\cC'$. Choose $D$ to be any
integer with $D > d_1 + \cdots + d_c - c + 1$. By
Remark \ref{socle} the minimal homogeneous generators of every ideal
in $\cC$ have degree less than $D$. Now let $(I',F')$ be an ideal
in $\cC'$ and let $C$ be a homogeneous complete intersection
of smallest type inside $(I',F')$, of height $c+1$. Since
$D$ exceeds the generator degrees of $I'$,
we can write $C = (C',G)$, where $C'$ is a homogeneous
complete intersection in $I'$ whose generators have minimal degree,
namely $d_1,\ldots,d_c$, and $G$ is a form in $(I', F')$ with $\deg G =
\deg F'=D$. Note that $G \notin I'$, so $(I',F') = (I',G)$. Let $J=C':I'$ be
the residual of $I'$ with respect to $C'$. Then $J \in \cC$ and  $C : (I', F') =
(C',G) : (I',G) = (J,G) \in \cC'$, as claimed.
\end{proof}

\medskip

The proof of Theorem \ref{two bdl} shows that the ideals considered
there generate a class $\cC$ satisfying the assumptions of
Proposition \ref{hypersurf sect}. Thus repeated application of the
proposition immediately gives the following consequence:

\smallskip

\begin{cor}
Let $k$ be an infinite field.  Then there are subschemes of any
codimension $c \geq 3$ in any
projective space $\mathbb P^d_k$ with $d \geq c$ that are homogeneously
licci, but not minimally
homogeneously licci.
\end{cor}

\medskip

We do not know if the result of Proposition \ref{hypersurf sect} is
still true if we
remove the condition that $\deg F \gg 0$.

\smallskip

\begin{remark}
 From the proof of Proposition \ref{hypersurf sect} we can be more
precise about what we mean by the hypothesis $\deg F \gg 0$.
Indeed, the numbers $d_1,\ldots,d_c$ cannot exceed the
Castelnuovo-Mumford regularity of $I$.  Hence it is certainly
enough to assume $\deg F \geq (\hbox{codim } I) \cdot (\hbox{reg
}I)$.
\end{remark}

\bigskip
\bigskip

\bibliographystyle{amsalpha}

\bigskip
\bigskip

\end{document}